\documentclass[12pt]{article}
\usepackage{color}
\usepackage{amsmath, amsthm, amsfonts, amssymb,latexsym}
\usepackage{graphicx} 
\usepackage{epstopdf}
\usepackage[bookmarks=false]{hyperref}
\usepackage{mathrsfs}
\UseRawInputEncoding

\newtheorem{thm}{Theorem}[section]
\newtheorem{cor}[thm]{Corollary}
\newtheorem{lem}[thm]{Lemma}

\theoremstyle{definition}
\newtheorem{defn}[thm]{Definition}

\newtheorem{exam}[thm]{Example}

\newtheorem{rem}[thm]{Remark}

\begin{document}

\def\prof{{\sc Proof.\ \ }}
\def\sect#1{\begin{center}\section{#1}\end{center}} 
\def\R#1{{\bf R}^{#1}}
\def\I#1#2{\int\limits_{#1}#2} 
\def\p#1{#1^{\prime}} 
\def\l#1#2#3{\lim_{#1 \rightarrow #2}#3}
\def\ld#1#2#3{\liminf_{#1 \rightarrow #2}#3} 
\def\lu#1#2#3{\limsup_{#1 \rightarrow #2}#3}
\def\E#1#2{{\bf E}^{#1}\left(#2\right)} 
\def\P#1#2{{\bf P}^{#1}\left(#2\right)} 
\def\Pt#1#2#3{{\bf P}(#1,#2,#3)} 
\def\pb#1{{\bf P}\left(#1\right)}
\def\eb#1{{\bf E}\left(#1\right)}
\def\Di#1{\,Dim\,#1}
\def\di#1{\,dim\,#1}
\def\fr#1#2{\frac{#1}{#2}}
\def\beq{\begin{equation}}
\def\eeq{\end{equation}}
\def\bea{\begin{eqnarray}}
\def\bean{\begin{eqnarray*}}
\def\eean{\end{eqnarray*}}
\def\eea{\end{eqnarray}}
\def\heq#1#2#3{\hbox to \hsize{\hskip #1 $#2$ \hss (#3)}}
\def\df#1#2{\frac{\displaystyle #1}{\displaystyle #2}}

\def\bdes{\begin{description}}
\def\ndes{\end{description}}

\newcommand{\bh}{{\bf h}}
\newcommand{\hf}{{\bf f}}
\newcommand{\he}{{\bf e}}
\newcommand{\hL}{{\bf L}}
\newcommand{\hg}{{\bf g}}
\newcommand{\hG}{{\bf G}}
\newcommand{\hM}{{\bf M}}

\newcommand{\bb}{\mathbb{b}}
\newcommand{\ww}{\mathbb{W}}
\newcommand{\hh}{\mathbb{H}}
\newcommand{\dd}{\mathbb{D}}
\newcommand{\cc}{\mathbb{C}}
\newcommand{\ee}{\mathbb{E}}
\newcommand{\zz}{\mathbb{Z}}
\newcommand{\nn}{\mathbb{N}}
\newcommand{\pp}{\mathbb{P}}
\newcommand{\qq}{\mathbb{Q}}
\newcommand{\ttt}{\mathbb{T}}

\def\rr{\Bbb R}
\def\a{\alpha}
\def\b{\beta}
\def\g{\gamma}
\def\s{\sigma}
\def\ep{\epsilon}
\def\d{\delta}
\def\D{\Delta}
\def\o{\omega}
\def\O{\Omega}
\def\n={\not=}
\def\u>{\wedge}
\def\d>{\vee}
\def\.{\bullet}
\def\l{\lambda}
\def\L{\Lambda}
\def\r{\rho}
\def\vf{\varphi}
\def\f{\phi}
\def\t{\tau}
\def\z{\zeta}
\def\na{\nabla}

\def\B{{I\!\!B}}
\def\E{{I\!\!E}}
\def\N{{I\!\!N}}
\def\P{{I\!\!P}}
\def\Q{{I\!\!\!Q}}
\def\R{{I\!\!R}}

\def\AA{\mathcal A}
\def\BB{\mathcal B}
\def\CC{\mathcal C}
\def\DD{\mathcal D}
\def\FF{\mathcal F}
\def\EE{\mathcal E}
\def\GG{\mathcal G}
\def\JJ{\mathcal J}
\def\LL{\mathcal L}
\def\NN{\mathcal N}
\def\PP{\mathcal P}
\def\SS{\mathcal S}

\def\Fht{\hbox{${\cal F}_{t}$}}
\def\Ght{\hbox{${\cal G}_{t}$}}

\def\U{\bigcup}
\def\Uu{\bigcap}
\def\Au{\forall}
\def\Eu{\exists}
\def\8u{\infty}
\def\0{\emptyset}
\def\mp{\longmapsto}
\def\rt{\rightarrow}
\def\up{\uparrow}
\def\dn{\downarrow}
\def\sub{\subset}
\def\vep{\varepsilon}

\def\ln{\langle}
\def\rn{\rangle}

\def\<<{\langle\!\langle}
\def\>>{\rangle\!\rangle}
\def\3|{|\!|\!|}

\def\ess{\text{\rm{ess}}}
\def\beg{\begin}

\def\beqt{\begin{equation}}
\def\neqt{\end{equation}}
\def\beq{\begin{equation}}
\def\neq{\end{equation}}

\def\bdes{\begin{description}}
\def\ndes{\end{description}}

\def\Ric{\text{\rm{Ric}}}
\def\Hess{\text{\rm{Hess}}}
\def\i{\text{\rm{i}}}
\def\ii{\text{\rm{ii}}}
\def\iii{\text{\rm{iii}}}
\def\iv{\text{\rm{iv}}}
\def\v{\text{\rm{v}}}
\def\vi{\text{\rm{vi}}}
\def\vii{\text{\rm{vii}}}
\def\viii{\text{\rm{viii}}}
\def\e{\text{\rm{e}}}

\def\Ra{\Rightarrow}
\def\Lra{\Leftrightarrow}
\def\rto{\longrightarrow}
\def\la{\leftarrow}
\def\ra{\rightarrow}
\def\ua{\uparrow}
\def\da{\downarrow}


\title{The optimal partition for multiparametric semialgebraic optimization
\thanks{Supported by the National Natural Science Foundation of China (11871118,12271061).}}

\author{Zi-zong Yan \thanks{School  of Information and Mathematics,
Yangtze University, Jingzhou, Hubei,
China(zzyan@yangtzeu.edu.cn).}
\ Xiang-jun Li \thanks{School  of Information and Mathematics,
Yangtze University, Jingzhou, Hubei,
China(franklxj001@163.com).}
\ and Jin-hai Guo\thanks{School  of Information and Mathematics,
Yangtze University, Jingzhou, Hubei,
China(xin3fei@21cn.com).}}
\date{}
\maketitle

\begin{abstract} In this paper we investigate the optimal partition approach for multiparametric conic linear optimization (mpCLO) problems in which the objective function depends linearly on vectors. We first establish more useful properties of the set-valued mappings early given by us \cite{YLG20} for mpCLOs, including continuity, monotonicity and semialgebraic property. These properties characterize the notions of so-called linearity and nonlinearity subsets of a feasible set, which serve as stability regions of the partition of a conic (linear inequality) representable set. We then use the arguments from algebraic geometry to show that a semialgebraic conic representable set can be decomposed into a union of finite linearity and/or nonlinearity subsets. As an application, we investigate the boundary structure of the feasible set of generic semialgebraic mpCLOs and obtain several nice structural results in this direction, especially for the spectrahedon.


\end{abstract}

\maketitle

\section{Introduction}
This paper is mainly concerned with the optimal partition approach for multiparametric conic linear optimization (mpCLO) problems by means of the connections existing between the vectors of parameters. In particular, we focus on the sensitivity analysis of semialgebraic mpCLOs.

Formally we consider the following pair of mpCLOs
\begin{equation} \label{primaljihe4} \begin{array}{ll} \min\limits_{x} & \langle c, x \rangle \\ s.t. & Ax=b, \\ & Mx=Md+v,
 \\ & x\in K \end{array} \end{equation} and
\begin{equation} \label{primaljihe3} \begin{array}{ll} \max\limits_{w,y} & b^Tw \\ s.t. & A^Tw+y=c+M^Tu, \\ &w\in\mathbb{R}^m, y \in K^*, \end{array} \end{equation} where $K\subset \mathbb{R}^q$ is a pointed, closed, convex, solid (with non-empty interior) cone (Ref. \cite{NN94}), and $K^*$ is the dual of $K$, that is, \[ K^*=\{y\in \mathbb{R}^q| \langle y,x\rangle\geq 0, \forall x\in K\}. \] Hence $d,c\in\mathbb{R}^q$, $b\in\mathbb{R}^m$ and $A \in \mathbb{R}^{m\times q}$ and $M\in \mathbb{R}^{r\times q}$ are fixed data, and $u,v\in\mathbb{R}^r$ are vectors of parameters. Such a pair of problems covers linear programming (LP) and semidefinite programming (SDP) problems. As an active research subject today, SDPs have gained a lot of attention from many practical applications in fields such as synthesis of filter and antennae arrays, structural design, stability analysis of mechanics, robust optimization, and relaxation of combinatorial optimization---topics of contemporary interest (Refs. \cite{BN01,BV04,LV98,SY15,Tod01}).

A astonishing aspect is that the Lagrangian dual of the problem (\ref{primaljihe3}) is the same as a lifting of the problem (\ref{primaljihe4}), and they can be expressed in one unified form:
\begin{equation} \label{primaljihe1} \begin{array}{ll} \min\limits_{x} & \langle c+M^Tu,x \rangle \\ s.t. & Ax=b,\\ & x \in K, \end{array} \end{equation} see \cite{YLG20}. In the same way, both the Lagrangian dual of the problem (\ref{primaljihe4}) and a lifting of the problem (\ref{primaljihe3}) can be described as another unified form:
\begin{equation} \label{primaljihe2} \begin{array}{ll} \max\limits_{w,s,y} & b^Tw + (Md+v)^Ts\\ s.t. & A^T w +M^Ts+y=c, \\ & w\in\mathbb{R}^m, s \in\mathbb{R}^r, y\in K^*. \end{array} \end{equation} Their optimality characterizations lead to {\it a multiparametric Karush-Kuhn-Tucker (mpKKT) property}. More precisely, a triple of vectors $(x;w,y)\in \mathbb{R}^{q\times r\times q}$ is optimal for the problems (\ref{primaljihe4}) and (\ref{primaljihe3}) if and only if it satisfies {\it the following system of conic linear
equations} \begin{subequations}
\begin{align}\label{kktcond1} &A x = b, \quad Mx =Md+v, \quad x \in K, \\ \label{kktcond2} &
A^Tw +y =c+M^Tu,\quad w\in\mathbb{R}^m, \quad y \in K^*,\\ \label{kktcond3} & \langle x,y \rangle=0,
\end{align} \end{subequations} (see Remark \ref{remkkt1} and Theorem \ref{classical1} below), in which the equality (\ref{kktcond3}) is well known {\it the complementary slackness property.} The surprising aspect of this result is that it holds for some vector pairs $(u,v)\in\mathbb{R}^{r\times r}$, under a Slater-type condition.

The mpKKT property (\ref{kktcond1})-(\ref{kktcond3}) is of interest because it covers the classical KKT property of either the primal-dual problems (\ref{primaljihe3}) and (\ref{primaljihe1}) or (\ref{primaljihe4}) and (\ref{primaljihe2}) for some vector pairs $(u,v)=(-s,v)$ (see also Remark \ref{remkkt1} below). This characterization result motivates us to explore in detail the relationship between the vectors of parameters $u$ and $v$. This relationship can be read as follows.

\defn \label{xinzd3} The set-valued mappings $\Phi,\Psi:\mathbb{R}^r \rightrightarrows \mathbb{R}^r$ are defined by
\begin{subequations}
\begin{align} \label{biaoshi1} \Phi(u)&=\{v\in\mathbb{R}^r| \exists \ (\bar{x};\bar{w},\bar{y})\ \operatorname{s.t. \ the\ system} \ (\ref{kktcond1})-(\ref{kktcond3})\ \operatorname{holds}\}, \\ \label{biaoshi2} \Psi(v)& = \{u \in\mathbb{R}^r | \exists \ (\bar{x};\bar{w},\bar{y})\ \operatorname{s.t. \ the\ system} \ (\ref{kktcond1})-(\ref{kktcond3})\ \operatorname{holds}\},
\end{align} \end{subequations} respectively. \upshape

We will show that under a Slater-type condition, the above definition are essentially equivalent to the one in the paper \cite{YLG20} (see Subsection \ref{oldnew1}). The latter was used by us to discuss the behaviour of the optimal partition of a conic representable set.

Clearly, the domains of the set-valued mappings $\Psi$ and $\Phi$ are contained in the following two sets
\begin{subequations}
\begin{align} \label{biaoshi10} \varTheta_P &=\{v\in\mathbb{R}^r| \exists \ \bar{x} \ \operatorname{ s.t. \ the \ system \ (\ref{kktcond1}) \ holds}\}, \\ \label{biaoshi20} \varTheta_D&=\{u\in\mathbb{R}^r| \exists \ (\bar{w},\bar{y}) \ \operatorname{ s.t. \ the \ system \ (\ref{kktcond2}) \ holds}\},
\end{align} \end{subequations} respectively, in which every element is naturally associated with the feasible sets of problems (\ref{primaljihe4}) and (\ref{primaljihe3}). One of the main reasons for considering them together is that these two sets enjoy more powerful results in the context of the optimal partition approach for mpCLOs (see Section \ref{parti14}).

Such sets can be represented as conic representable sets, including polyhedra and spectrahedra (see Lemma \ref{conicr1} and the corresponding remark). Polyhedra enjoys considerable interest throughout pure and applied mathematics, and spectrahedra inherit many of the favorable properties of polyhedra. From the viewpoint of geometry, the boundary structure of spectrahedra is far more intricate than that of polyhedra. Recently there has been considerable interest in understanding the geometry of spectrahedra (Refs. \cite{GW15,HV07,Hen11,Las17,NP10,Nie08,San15,Sch18}). Unfortunately, except for polyhedra, relatively little is known about the boundary structure of spectrahedra.

\begin{figure}[htbp] \centering\includegraphics[width=4.0in]{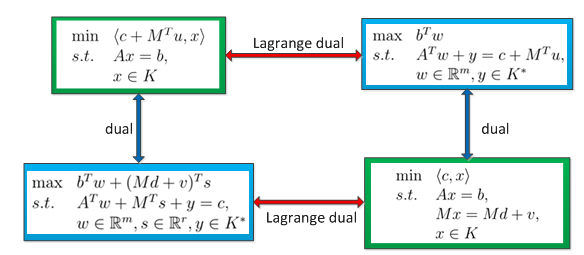} \caption{Various different duals}\label{fig:1} \end{figure}

Like the classical KKT property, the mpKKT property (\ref{kktcond1})-(\ref{kktcond3}) characterizes the optimality characteristics of the problems (\ref{primaljihe4}) and (\ref{primaljihe3}) whenever the complement slackness property holds (see Theorem \ref{classical1} below). So in our terminology, the problems (\ref{primaljihe4}) and (\ref{primaljihe3}) are called {\it a pair of primal-dual mpCLOs} (actually, their lifting models (\ref{primaljihe2}) and (\ref{primaljihe1}) is another pair of primal-dual mpCLOs, see Figure \ref{fig:1}). Correspondingly, the sets $\varTheta_{P}$ and $\varTheta_{D}$ are called {\it a pair of primal-dual conic representable sets}.

The following is assumed to hold throughout this paper.

{\bf Assumption 1}. Both the matrices $A$ and $M$ are of full rank, and the range spaces of their transpose don't intersect, i.e., $R(A^T)\cap R(M^T)=\emptyset$.

{\bf Assumption 2}. There exists a (primal) feasible point $x\in \operatorname{int}(K)$ for the problem (\ref{primaljihe1}), and there exists a (dual) feasible point $(w,s,y)$ with $y\in \operatorname{int}(K^*)$ for the problem (\ref{primaljihe2}).

Assumption 2 implies (see e.g. \cite{YLG20}) that both $\varTheta_{P}$ and $\varTheta_{D}$ satisfy a Slater type condition, i.e., they have at least an interior point in $\mathbb{R}^r$.

It should be noted that column vectors of $M^T$ denote perturbation directions of the objective function of (\ref{primaljihe1}). On the other hand, we \cite{YLG20} assumed that $d$ and $(0,0,c)$ are feasible for the problems (\ref{primaljihe1}) and (\ref{primaljihe2}), respectively. Although this assumption is unnecessary for this paper, for convenience we still adopt the same notation of the paper \cite{YLG20} to preserve the vector $d$.

Motivated by our recent work \cite{YLG20}, the main subject of this paper is to make further study the optimal partition of the primal and dual conic representable sets defined by (\ref{biaoshi10})-(\ref{biaoshi20}). For the optimal partition, an overview of the literature reveals that there are two major types of invariancy sets---{\it linearity and nonlinearity sets}---whose difference is whether the optimal objective value functions of the problems (\ref{primaljihe4}) and (\ref{primaljihe3}) are linear with respect to the vector of parameter $u$ or $v$ (see Theorems \ref{linritys1} and \ref{linritys2}).

There are significant differences between polyhedra and specthedra although specthedra naturally generalize the class of polyhedra. Recall that an extreme point of a convex set is called a vertex if its normal cone is full-dimensional, and is called smooth if its normal cone is one-dimensional. For a polyhedron, extreme points and vertices coincide, and there are only finitely many of them. However, the unit ball $\{x\in\mathbb{R}^q|\|x\|\leq 1\}$, which is affinely isomorphic to a specthedron, has uncountably many extreme points and no vertices when $q\geq 2$ (Ref. \cite{MST16}). Linearity and nonlinearity sets facilitate the study of their boundary structure. For example, a linearity set is usually used to analyze shadow prices or marginal revenues for the sensitivity analysis of LPs (Refs. \cite{GS55,Gal86}), whereas a nonlinearity set is in general regarded as a stability region and its identification has a great influence on the post-optimal analysis of SDPs (Refs. \cite{CA17,CA18}).

For theoretical and practical significance the question of the distinction between linearity and nonlinearity sets in the description of a conic representable set is very important. Although the general question has so far been elusive, many results have been obtained in answer of the question for some special cases. For example, a linearity set is known to be convex for polyhedra \cite{GG08}, and a nonlinearity set is open for spectrahedra \cite{YLG20}; a polyhedron has only linearity subsets, and a specthedron usually has a nonlinearity subset (see Corollary \ref{ghui1} and the corresponding remarks). Such two types of invariancy sets were first distinguished by Mohammad-Nezhad and Terlaky \cite{MT20} in the context of SDPs, although the first type has been studied extensively in LPs (Refs. \cite{AD92,Deh07,Gre94,JR03}), quadratic programming and linear complementarity problems (Ref. \cite{Ber97}).

In addition, multiparametric programming has emerged in the last two decades as an important optimization-based tool for systematically analyzing the effect of uncertainty and variability in optimization. There are plenty of works on multiparametric linear (Refs. \cite{BK21,Gal95,GT72,Fil97,Sch87}), nonlinear (Refs.\cite{Fia76,Fia83,FK86,Koj79}), quadratic (Refs. \cite{BMP02,BD95,PD02,WW20}), mixed-integer linear programming (Refs. \cite{DP00,GN77,JM06,MM75}) and others (Refs. \cite{CJ91,MC80}). To the best of our knowledge, in the general case there are hardly any results regarding the question mentioned above; and even for the case of the two parameters ($r=2$), the answer is unknown.

In this paper we exhibit more algebraic and geometric properties of the linearity and the nonlinearity sets than in the literature. They are based on supposedly good properties of the set-valued mappings $\Phi$ and $\Psi$. We show that if $K$ is a semialgebraic set, then both $\Phi$ and $\Psi$ are semialgebraic mappings. That is, their graphs $gph(\Phi)$ and $gph(\Psi)$ on $\mathbb{R}^r\times \mathbb{R}^r$ are semialgebraic---those whose graphs can be written as a finite union of sets, each obtained by finitely many polynomial equalities and inequalities. We then use the arguments from algebraic geometry, in particular properties of semialgebraic sets, to show that a semialgebraic conic representable set can be decomposed into a union of finite linearity and/or nonlinearity subsets. Our results imply that the feasible set has finitely vertices, and in particular the spectrahedon usually has a smooth surface. These geometric properties have the practical significance in describing the boundary structure of the feasible set, especially for the spectrahedon. All these lead to a transparent and unified treatment of the boundary structure of the feasible sets of semialgebraic optimization problems, covering, in particular, polynomial optimization problems, LPs and SDPs.

The organization of this paper is as follows. In Section 2, we recall some preliminary results from set-valued analysis and semialgebraic geometry. In Section 3, we review the essential facts from \cite{YLG20} about the set-valued mapping and discuss their continuity and monotonicity, and prove that they are semialgebraic mappings when $K$ is a semialgebraic set. In Section 4, we investigate the sensitivity of the optimal partition of semialgebraic conic representable sets and obtain several nice structural results in this direction. The main tool of our investigation comes from semialgebraic geometry. The last section contains some conclusions.

\section{Preliminaries}
\subsection{Basic notation, terminology and facts}
Throughout this paper we follow the terminology and notation of the book \cite{SPR06,Roc70,RW09} and the paper \cite{YLG20}.

The symbol $\mathbb{R}^q$ denotes a $q$-dimensional Eucildean space with the standard inner product $\langle \cdot,\cdot\rangle$ and the corresponding norm $\|x\|$. A subset $C\subset \mathbb{R}^q$ is a convex if for any $x,y\in C$ and any $\lambda\in [0,1]$, $\lambda x+(1-\lambda)y\in C$, and a point $z\in C$ is an extreme point of $C$ if there is no any $\lambda\in(0,1)$ and $x,y\in C$ such that $z=\lambda x+(1-\lambda)y$.

A subset $D\subset \mathbb{R}^q$ is a cone if it is closed under positive scalar multiplication. It is called pointed if $D$ contains no straight line, i.e., $x\in D, x\ne0$ imply $-x\notin D$.

For any subset $W\subset \mathbb{R}^q$, denote $\operatorname{int}(W)$, $\operatorname{cl}(W)$ and $\operatorname{conv}(W)$ as the interior, closure and convex hull of $W$, respectively.

Let $C$ be a nonempty convex set in $\mathbb{R}^q$. Given a boundary point $\bar{x}$ of $C$, its normal cone $\operatorname{Normal}(C,\bar{x})$ is defined by
\[ \operatorname{Normal}(C,\bar{x})=\{c\in \mathbb{R}^q | \langle c,\bar{x}\rangle\geq \langle c,x\rangle \ \operatorname{for \ all} \ x\in C\}.\]

For any $A\in\mathbb{R}^{m\times q}$, the notation $\operatorname{rank}(A)$ and $R(A)$ denote the rank and the range space of the matrix $A$. In particular, if $m=q$ and $A^2=A=A^T$, then $A$ is called a symmetric projection matrix.

Note that the feasible sets of the problems (\ref{primaljihe1}) and (\ref{primaljihe2}), denoted by $\mathscr{X}$ and $\mathscr{Y}$, respectively, do not depend on the given vectors of parameters $u$ and $v$. Without causing confusion, the optimal solutions are denoted by $x^*(u)$ and $(w^*(v),s^*(v),y^*(v))$, respectively, if they exist. That is,
\begin{eqnarray*} x^*(u) &\in& \mathscr{X}^*(u)= \arg\min \{\langle c+M^Tu,x\rangle | x\in \mathscr{X}\}, \\ (w^*(v),s^*(v),y^*(v)) & \in & \mathscr{Y}^*(v)= \arg\max \{b^Tw + (Md+v)^Ts |(w,s,y)\in \mathscr{Y}\}. \end{eqnarray*} For brevity, sometimes, we also replace $ (w^*(v),s^*(v),y^*(v)) \in \mathscr{Y}^*(v)$ by $y^*(v) \in \mathscr{Y}^*(v)$.

A conic representable set is the solution of a conic linear inequality. The following lemma shows that both $\varTheta_P$ and $\varTheta_D$ are conic representable sets, in which $q-m-r=0$ is allowed.

\lem \label{conicr1} (1) $\varTheta_D$ can be rewritten as
\begin{equation} \label{cong1} \varTheta_D = \{ u\in\mathbb{R}^r | \exists w\in\mathbb{R}^m \ \operatorname{s.t.}\ c+M^Tu-A^Tw\in K^* \}. \end{equation}

(2) Suppose that three range spaces $R(A^T)$,$R(M^T)$, $R(B^T)$ are orthogonal to each other, where $B\in\mathbb{R}^{(q-m-r)\times q}$. If $d\in \{x\in\mathbb{R}^q| Ax=b\}$, then $\varTheta_P$ can be redescribed by
 \begin{equation} \label{cong2} \varTheta_P = \{v\in\mathbb{R}^r | \exists t\in\mathbb{R}^{q-m-r} \ \operatorname{s.t.}\ d+M^T(MM^T)^{-1} v-B^Tt\in K \}. \end{equation} \upshape

\begin{proof} The second claim references \cite[Corollary 2.4]{YLG20}, while the first claim follows from (\ref{kktcond2}) trivially. \end{proof}

When $K$ is the nonnegative orthant, the sets described as (\ref{cong1}) and (\ref{cong2}) are linear representable sets or polyhedra. When $K$ is the cone of symmetric positive semidefinite matrices, they are linear matrix representable sets or (projected) spectrahedra (Refs. \cite{AL12,BN01,BPT13,Nem07}). Such two types of sets represent the geometry underlying LPs and SDPs, respectively.

\subsection{Variational analysis}
Let $f:\mathbb{R}^q\rightarrow \mathbb{R}^l$ denote a mapping from $\mathbb{R}^q$ to $\mathbb{R}^l$. The preimage of a subset $T$ of $\mathbb{R}^l$ under the mapping $f$ is defined by
\[f^{-1}(T) := \{x\in \mathbb{R}^q | f(x)\in T \}. \]

A set-valued mapping $F$ from $\mathbb{R}^q$ to $\mathbb{R}^l$, denoted by $F:\mathbb{R}^q\rightrightarrows \mathbb{R}^l$, is a mapping
from $\mathbb{R}^q$ to the power set of $\mathbb{R}^l$. Thus for each point $x\in\mathbb{R}^q$, $F(x)$ is a subset of $\mathbb{R}^l$.
For a set-valued mapping $F:\mathbb{R}^q\rightrightarrows \mathbb{R}^l$, the domain, graph and inverse mapping of $F$ are defined by
\begin{eqnarray*} dom F &:=& \{x \in \mathbb{R}^q | F(x) \ne \emptyset \},\\
gph F &:=& \{(x, z) \in \mathbb{R}^q\times \mathbb{R}^l | z\in F(x)\}, \\ F^{-1}(z)&:=&\{x\in\mathbb{R}^q| z\in F(x)\}\end{eqnarray*}
and the image of $F$ over the subset $D$ of $\mathbb{R}^q$ is described as
\[F(D):=\{F(x)|x\in D\}. \]

A subset $G$ of $\mathbb{R}^q\times \mathbb{R}^q$ is called monotone if for each choice of $(x_1, z_1)$ and $(x_2,z_2)$ in $G$ one has $\langle x_1-x_2,z_1-z_2\rangle\geq 0$. A set-valued mapping from $\mathbb{R}^q$ to $\mathbb{R}^q$ is monotone if its graph is a monotone set. For more on monotonic set-valued mappings see e.g.\cite[Chapter II]{Bre73}.

In the following we briefly review the continuity of the set-valued mappings from \cite{Roc14,RW09}.

{\bf(Painlev$\acute{e}$-Kuratowski limits of sequences)} We first recall basic notions about (Painlev$\acute{e}$-Kuratowski) limits of sets.
Given a sequence $\{C_{\nu}\}_{\nu\in\mathbb{N}}$ of subsets of $\mathbb{R}^q$ we define:

$\bullet$ the outer limit $\limsup\limits_{\nu\rightarrow\infty}C_{\nu}$, as the set of all accumulation points of sequences $\{x_{\nu}\}_{\nu\in\mathbb{N}}\subset \mathbb{R}^q$
with $x_{\nu}\in C_{\nu}$ for all $\nu\in\mathbb{N}$. In other words, $x\in\limsup\limits_{\nu\rightarrow\infty}C_{\nu}$ if and only if for every neighbourhood $U$ of $x$ and
$N\geq 1$ there exists $\nu\geq N$ with $C_{\nu}\cap U\ne \emptyset$;

$\bullet$ the inner limit $\liminf\limits_{\nu\rightarrow\infty}C_{\nu}$,, as the set of all limits of sequences $\{x_{\nu}\}_{\nu\in\mathbb{N}}\subset \mathbb{R}^q$ with $x_{\nu}\in C_{\nu}$ for all $\nu\in\mathbb{N}$. In other words, $x\in\liminf\limits_{\nu\rightarrow\infty}C_{\nu}$ if and only if for every neighbourhood $U$ of $x$ there exists $N\geq 1$ such that for all $\nu\geq N$ we have $C_{\nu}\cap U \ne \emptyset$.

Furthermore, we say that the limit of the sequence $\{C_{\nu}\}_{\nu\in\mathbb{N}}$ exists if the outer and inner limit sets are equal. In this case we write
\[\lim\limits_{\nu\rightarrow\infty}C_{\nu}=\limsup\limits_{\nu\rightarrow\infty}C_{\nu}=\liminf\limits_{\nu\rightarrow\infty}C_{\nu}. \]

\defn A set-valued mapping $F:\mathbb{R}^q\rightrightarrows \mathbb{R}^l$ is called outer semicontinuous at $\bar{x}$ if
\[ \limsup\limits_{x\rightarrow \bar{x}} F(x)\subset F(\bar{x}), \] or equivalently, $\limsup\limits_{x\rightarrow \bar{x}} F(x) = F(\bar{x})$, and inner semicontinuous at $\bar{x}$ if
\[ \liminf\limits_{x\rightarrow \bar{x}} F(x)\supset F(\bar{x}), \] or equivalently, $\liminf\limits_{x\rightarrow \bar{x}} F(x) = F(\bar{x})$. It is called continuous at $\bar{x}$ if both conditions hold, i.e., if $F(x)\rightarrow F(\bar{x})$ as $x\rightarrow \bar{x}$. \upshape


Notice that every outer semicontinuous set-valued mapping has closed values. In particular, it is well known that

$\bullet$ $F:\mathbb{R}^q\rightrightarrows \mathbb{R}^l$ is outer semicontinuous if and only if it has a closed graph.

When $f:\mathbb{R}^q\rightarrow \mathbb{R}^l$ is a single-valued function, both outer and inner semicontinuity reduce to the standard notion of continuity.

\subsection{Semialgebraic geometry}
A subset $S$ of $\mathbb{R}^q$ is called semialgebraic if it is a finite union of sets of the form
\[ \{x\in \mathbb{R}^q | f_i(x)=0, i=1,2,\cdots,k; f_i(x)>0, i=k+1,\cdots,l\}, \] where all $f_i$
are polynomials. In other words, $S$ is a union of finitely many sets, each defined by finitely many polynomial equalities and inequalities. The following simple facts will be used.

$\bullet$ Both the nonnegative orthant and the cone of symmetric positive semidefinite matrices are semialgebraic.

$\bullet$ If $S\subset\mathbb{R}^q$ is semialgebraic, then so are its interior and closure.

$\bullet$ If $S\subset\mathbb{R}^q$ and $T\in\mathbb{R}^l$ are semialgebraic, then their Cartesian product $S\times T=\{(x,z)|x\in S, z\in T\}$ is a semialgebraic set of $\mathbb{R}^q\times \mathbb{R}^l$.

A semialgebraic set $S\subset\mathbb{R}^q$ is semialgebraically connected if $S$ is not the disjoint union of two non-empty semialgebraic sets that are both closed
in $S$. Or equivalently, $S$ does not contain a non-empty semialgebraic strict subset which is both open and closed in $S$.

A (semialgebraic) set $S\subset\mathbb{R}^q$ is (semialgebraically) path connected when for every $x, y$ in $S$, there exists a (semialgebraic) path from $x$ to $y$,
i.e., a continuous (semialgebraic) function $\phi: [0, 1] \rightarrow S$ such that $\phi(0) = x$ and $\phi(1) = y$. Sometimes, the path connected is called simply connected.

The following three results reference \cite[Theorems 3.7, 5.22 and 5.23]{SPR06}.

\thm \label{conec1} If $S\in\mathbb{R}^q$ is semialgebraic and convex then it is semialgebraic connected. \upshape

\thm \label{conec2} A semialgebraic set $S$ of $\mathbb{R}^q$ is semialgebraically connected if and only if it is connected. Every semialgebraic set (and in particular every semialgebraic subset) of $\mathbb{R}^q$ has a finite number of connected components, each of which is semialgebraic. \upshape

\thm \label{conec3} A semialgebraic set is semialgebraically connected if and only if it is semialgebraically path connected. \upshape

A mapping $F:\mathbb{R}^q\rightrightarrows \mathbb{R}^l$ is said to be semialgebraic if its graph is a semialgebraic subset in $\mathbb{R}^q\times \mathbb{R}^l$. A fundamental fact about semialgebraic sets is provided by the Tarski-Seidenberg Theorem \cite[Theorem 2.3]{Cos02}.

\thm \label{TST1} {\bf (Tarski-Seidenberg theorem)} The image of a semialgebraic set by a semialgebraic mapping is semialgebraic. \upshape



A mapping $f:\mathbb{R}^q \rightarrow \mathbb{R}^q$ is a semialgebraically homeomorphism if it is a homeomorphism of $\mathbb{R}^q$ onto $\mathbb{R}^q$ and both $f$ and $f^{-1}$ are semialgebraically, respectively. In particular, for $f$ to be a semialgebraically homeomorphism it does not suffice for it to be a homeomorphism that is semialgebraically: the inverse mapping must also be semialgebraically. Henceforth, we use the symbol $\cong$ to indicate that two semialgebraic sets are semialgebraically homeomorphic.

The following Hardt's triviality theorem appeared in \cite{har80}.

\thm \label{hards1} {\bf (Hardt's triviality theorem)} Let $S\in\mathbb{R}^q$ and $T\in \mathbb{R}^l$ be semialgebraic sets. Given a continuous semialgebraic mapping $f:S\rightarrow T$, there exists a finite partition of $T$ into semialgebraic sets $T=\bigcup\limits_{i=1}^kT_i$, so that for each $i$ and any $x_i\in T_i$, \[T_i\times f^{-1}(x_i)\cong f^{-1}(T_i).\] \upshape

By applying Hardt's triviality theorem to the projection mapping, we can get a semialgebraic partition of the semialgebraic set, which will provide a theoretical foundation for the optimal partition of conic representable sets.

\thm \label{TST2} Let $F:\mathbb{R}^q\rightrightarrows \mathbb{R}^l$ be a semialgebraic set-valued mapping. Then there exists a partition of the domain of $F$ into semialgebraic sets $X_1, X_2,\cdots, X_k$, so that for each $i$ and any $x_i\in X_i$, \[X_i\times (\{x_i\}\times F(x_i))\cong X_i\times F(X_i).\] \upshape

\begin{proof} We follow partly from the proof of \cite[Corollary 2.27]{DL13}. Let $p:gph F\rightarrow \mathbb{R}^q$ be the projection onto the first $q$ coordinates. By applying Theorem \ref{hards1} to $p$, we get a partition of the domain of $F$ into semialgebraic sets $X_1,X_2,\cdots,X_k$, so that for each $i$ and any $x_i\in X_i$, $X_i\times p^{-1}(x_i)=X_i\times (\{x_i\}\times F(x_i))$ is semialgebraically homeomorphic to $p^{-1}(X_i)=X_i\times F(X_i)$. The proof is finished. \end{proof}

\section{Review and extension: set-valued mappings} \label{review3}
The set-valued mappings between the primal and dual conic representable sets were first presented and analyzed in \cite{YLG20}. This section presents a brief review of the mappings and gets more important and useful properties, including continuity, monotonicity and semialgebraic property. These review and extension will be key ingredients for many results presented in the next section.

\subsection{New and old definitions} \label{oldnew1}
Let us start with a useful intermediate result that comes from Definition \ref{xinzd3}.

\lem \label{techical1} Let $(u,v)\in \mathbb{R}^r\times \mathbb{R}^r$ and $(\bar{x};\bar{w},\bar{s},\bar{y})\in\mathbb{R}^q\times \mathbb{R}^m\times \mathbb{R}^r\times \mathbb{R}^q$ be given. Then the following three statements are equivalent:

(1) $(\bar{x};\bar{w},\bar{y})$ is a solution of the system (\ref{kktcond1})-(\ref{kktcond3});

(2) $(\bar{x}; \bar{w},\bar{y})$ is a solution of the following system
 \begin{subequations}
\begin{align}\label{kktcond10} &A x = b, \quad x \in K, \\ \label{kktcond20} &
A^T w + y= c+M^Tu,\quad w \in\mathbb{R}^m, \quad y \in K^*,\\ \label{kktcond30} & \langle x, y\rangle=0
\end{align}
\end{subequations} and $v\in\Phi(u)$;

(3) $(\bar{x};\bar{w},\bar{s},\bar{y})$ is a solution of the following system
 \begin{subequations}
\begin{align}\label{fjkktcond10} &A x = b, \quad M x=Md+v, \quad x \in K, \\ \label{fjkktcond20} &
A^Tw + M^Ts+ y = c ,\quad w \in\mathbb{R}^m, \quad s \in\mathbb{R}^r, \quad y \in K^*,\\ \label{fjkktcond30} & \langle x,y \rangle=0
\end{align} \end{subequations} and $u=-\bar{s}\in\Psi(v)$. \upshape

\rem \label{remkkt1} It is known that the system (\ref{kktcond10})-(\ref{kktcond30}) denotes the KKT property of the primal-dual pair of problems (\ref{primaljihe1})-(\ref{primaljihe3}) (for example, see \cite{BN01,BPT13,BL05,BV04}). That is, $(\bar{x}; \bar{w},\bar{y})$ is a solution of the system (\ref{kktcond10})-(\ref{kktcond30}) if and only if it is a pair of optimal solutions such that \begin{equation}\label{setdef1} M(\bar{x}-d)\in \Phi(u).\end{equation} Analogously, the system (\ref{fjkktcond10})-(\ref{fjkktcond30}) denote the KKT property of the primal-dual pair of problems (\ref{primaljihe4})-(\ref{primaljihe2}); and $(\bar{x}; \bar{w}, \bar{s}, \bar{y})$ is a solution of the system (\ref{fjkktcond10})-(\ref{fjkktcond30}) if and only if it is a pair of optimal solutions such that \begin{equation}\label{setdef2} -\bar{s}\in \Psi(v). \end{equation}
So Lemma \ref{techical1} gives the following equivalence:
\begin{eqnarray*} \label{djgxi1}
 && \operatorname{The \ KKT \ property} \ (\ref{kktcond10})-(\ref{kktcond30}) \ \operatorname{holds}\ +\ v\in\Phi(u)\\ \label{djgxi2} &\Leftrightarrow& \operatorname{The \ KKT \ property} \ (\ref{fjkktcond10})-(\ref{fjkktcond30}) \ \operatorname{holds}\ +\ (u=-s)\in\Psi(v) \\ \label{djgxi3} & \Leftrightarrow& \operatorname{The \ system} \ (\ref{kktcond1})-(\ref{kktcond3}) \ \operatorname{holds.} \end{eqnarray*} \upshape

 This remark yields the optimality of the problems (\ref{primaljihe4})-(\ref{primaljihe3}) and (\ref{primaljihe1})-(\ref{primaljihe2}).

\thm \label{classical1} One of the three arguments of Lemma \ref{techical1} holds if and only if the following three conditions hold simultaneously:

 (1) Both the problems (\ref{primaljihe1}) and (\ref{primaljihe4}) have the same optimal solution $\bar{x}$;

 (2) Both the problems (\ref{primaljihe2}) and (\ref{primaljihe3}) have the same optimal slackness vector $\bar{y}$;

 (3) $(\bar{x},\bar{y})$ satisfies the complement slackness property (\ref{kktcond3}).
\upshape

From Theorem \ref{classical1}, we say that the system (\ref{kktcond1})-(\ref{kktcond3}) denotes the KKT property of the pair of primal-dual problems (\ref{primaljihe4})-(\ref{primaljihe3}) (or (\ref{primaljihe1})-(\ref{primaljihe2})).

\lem \label{invers1} $\Phi^{-1} = \Psi.$ That is, $u\in\Psi(v)$ if and only if $v\in\Phi(u)$. \upshape

\begin{proof} By Definition \ref{xinzd3}, $\Psi(v)$ is well defined and $u\in\Psi(v)$ if and only if $(\bar{x};\bar{w},\bar{y})$ is a solution of the system (\ref{kktcond1})-(\ref{kktcond3}), or equivalently, $\Psi(u)$ is well defined and $v\in\Phi(u)$. \end{proof}

From Lemma \ref{invers1}, the above equivalence can be redescribed as follows
 \begin{eqnarray*} \label{djgxi1}
 && \operatorname{The \ KKT \ properties} \ (\ref{kktcond10})-(\ref{kktcond30}) \ and \ (\ref{fjkktcond10})-(\ref{fjkktcond30}) \operatorname{with} u=-s \ \operatorname{hold}\\ &\Leftrightarrow& v\in\Phi(u) \\ &\Leftrightarrow& u\in\Psi(v) \\ \label{djgxi3} & \Leftrightarrow& \operatorname{The \ mpKKT \ property} \ (\ref{kktcond1})-(\ref{kktcond3}) \ \operatorname{holds.} \end{eqnarray*} So, under the set-valued mappings, the three KKT properties are essentially equivalent to each other.

\begin{figure}[htbp]
\centering\includegraphics[width=4.5in]{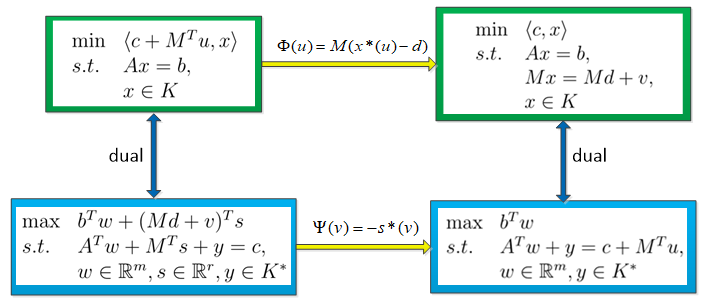}
\caption{Set-valued mappings}\label{fig:2}
\end{figure}

When the complement slackness property (\ref{kktcond3}) holds, from (\ref{setdef1}) and (\ref{setdef2}), the set-valued mappings $\Phi$ and $\Psi$ can be rewritten as follow
 \begin{subequations}
\begin{align}\label{ghjj1} \Phi(u) &= \{ M(x^*(u)-d)| x^*(u) \in \mathscr{X}^*(u) \}, \\ \label{ghjj2} \Psi(v) &= \{-s^*(v)| (w^*(v),s^*(v),y^*(v)) \in \mathscr{Y}^*(v) \},
\end{align} \end{subequations} respectively. Such a definition, first proposed by us \cite{YLG20}, establishes the connections between either the problems (\ref{primaljihe1}) and (\ref{primaljihe4}) by $\Phi$ or (\ref{primaljihe2}) and (\ref{primaljihe3}) by $\Psi$, see Figure \ref{fig:2}. It is easy to verify that Lemma \ref{techical1} and Theorem \ref{classical1} still hold for the definition.

By comparing (\ref{biaoshi1})-(\ref{biaoshi2}) and (\ref{ghjj1})-(\ref{ghjj2}), we can see the two definitions of the set-valued mappings are pretty close. The main gaps are whether the complement slackness property (\ref{kktcond3}) holds. In general, the old definition (defined by (\ref{ghjj1})-(\ref{ghjj2})) does not require this property. However, when the property (\ref{kktcond3}) holds, two definitions are consistent (see Remark \ref{remkkt1}). In the following statement, we always use the new definition (given by (\ref{biaoshi1})-(\ref{biaoshi2})) unless specifically mentioned.

\thm \label{slaterk1} Let $\Phi$ be defined by (\ref{ghjj1}). If $u\in \operatorname{int}(\varTheta_D)$, then $\Phi(u)$ is well defined, i.e., there is a parametric vector $v\in\varTheta_P$ such that $v\in \Phi(u)$. An analogous result holds for $\Psi$ defined by (\ref{ghjj2}). \upshape

A proof of Theorem \ref{slaterk1} is given by us \cite{YLG20} by using the strong duality theorem of conic programming. Indeed, if $u\in\operatorname{int}(\varTheta_D)$, then a pair of primal-dual problems (\ref{primaljihe1}) and (\ref{primaljihe3}) satisfy the Slater condition. By the strong duality theorem (for example, see \cite{BN01,BPT13,BL05,BV04}), there is $(\bar{x}; \bar{w},\bar{y})$ is a solution of the system
(\ref{kktcond10})-(\ref{kktcond30}); and if $v=M\bar{x}-Md$, then $v\in\Phi(u)$. By Remark \ref{remkkt1}, there is a triple of vectors $(\bar{x}; \bar{w},\bar{y})\in\mathbb{R}^q\times \mathbb{R}^r\times \mathbb{R}^q$ such that the mpKKT property (\ref{kktcond1})-(\ref{kktcond3}) holds. The same result applies to the set-valued mapping $\Psi$.

We give a new proof in Appendix A.

\rem \label{strong1} Theorem \ref{slaterk1} is essentially equivalent to the strong duality theorem of conic programming. By Remark \ref{remkkt1} and Theorem \ref{slaterk1}, if $u\in\operatorname{int}(\varTheta_D)$ or $v\in\operatorname{int}(\varTheta_P)$, then there is a triple of vectors $(\bar{x}; \bar{w},\bar{y})\in\mathbb{R}^q\times \mathbb{R}^r\times \mathbb{R}^q$ such that the mpKKT property (\ref{kktcond1})-(\ref{kktcond3}) holds, i,e, the strong duality theorem holds. Conversely, if the strong duality theorem holds, i.e., the mpKKT property (\ref{kktcond1})-(\ref{kktcond3}) holds, then from Remark \ref{remkkt1}, $\Phi$ and $\Psi$ given by (\ref{ghjj1}) and (\ref{ghjj2}) are well defined. \upshape

The Slater type condition in Theorem \ref{slaterk1} is necessary. We illustrate it with a semidefinite system, with $\mathbb{R}^{\frac{n(n+1)}{2}}\simeq\mathbb{S}^n$ the set of order $n$ symmetric matrices and $K=K^*=\mathbb{S}^n_+$ as the set of order $n$ symmetric positive semidefinite matrices. The inner product of $\mathbb{S}^n$ is $c \bullet d= \langle c,d\rangle=tr(cd)$. The following example is taken from the paper \cite{YLG20}.

\exam \label{patex10} Consider the single parameter SDP problem
\[ \begin{array}{ll} \min & x_{11}+ux_{22} \\ s.t. & x_{12} =1, \\ & \left(\begin{array} {cc} x_{11} & x_{12} \\ x_{12} & x_{22} \end{array} \right) \in\mathbb{S}^2_+, \end{array} \] in which $\varTheta_D=[0,+\infty]$. If $u>0$, then the problem has the optimal solution \[ x^*(u)=\left( \begin{array}{cc}\sqrt{u} &1\\ 1&\frac{1}{\sqrt{u}} \end{array}\right) \]with the dual optimal slackness matrix \[ \bar{y}^*(u)=\left( \begin{array}{cc}1 &-\sqrt{u}\\ -\sqrt{u}& u \end{array}\right), \] and the complement slackness property holds. Therefore, for any $u\in (0,+\infty)$,
\[\Phi(u) = \left(\begin{array} {cc} 0 & 0 \\ 0 & 1 \end{array} \right)\bullet \left( \begin{array}{cc}\sqrt{u} &1\\ 1&\frac{1}{\sqrt{u}} \end{array}\right)-\left(\begin{array} {cc} 0 & 0 \\ 0 & 1 \end{array} \right)\bullet \left( \begin{array}{cc}d_{11} &1\\ 1& d_{22}\end{array}\right)
=\frac{1}{\sqrt{u}}-d_{22} \] is well defined. However, for $u=0$, the problem is not solvable such that $\Phi(0)$ is undefined since $0\notin\operatorname{int}\varTheta_D$, i.e., the Slater type condition does not hold, although its Lagrange dual has 0 maximum at \[\bar{y}^*(0)=\left( \begin{array}{cc}1 &0\\ 0&0\end{array}\right)\] and there is no duality gap. \upshape

\rem \label{remgeo1} An interesting geometric interpretation of the set-valued mappings defined by (\ref{ghjj1})-(\ref{ghjj2}) was given by us \cite{YLG20}. For example, the set-valued mapping $\Phi(u)$ from $\varTheta_D$ to $\varTheta_P$ bind the problems (\ref{primaljihe1}) and (\ref{primaljihe4}) together, in which their perturbations occur in the objective function and in the right side hand, respectively. Firstly, $\varTheta_D$ denotes the objective perturbation set of the problem (\ref{primaljihe1}): geometrical, for a fixed $u\in\varTheta_D$, the hyperplane \begin{equation} \label{zcsup1} H_u=\{x\in \mathbb{R}^q| \langle c+M^Tu,x\rangle = \langle c+M^Tu,x^*(u)\rangle \}\end{equation} passes through the extreme point $x^*(u)$ and supports the feasible set $\mathscr{X}$. In particular, if $\Phi(u)$ is a set, then so is $H_u\cap \mathscr{X}$. Secondly, $\varTheta_P$ denotes the right side hand perturbation set of the problem (\ref{primaljihe4}): geometrical, the affine set \begin{equation} \label{zcsup2} C_v= \{x\in \mathbb{R}^q| Mx=Md+v\}\end{equation} passes through the extreme points $\bar{x}^*(v)\in H_u\cap \mathscr{X}$ and cuts the feasible set $\mathscr{X}$ for all $v\in\Phi(u)$. Finally, the optimal objective function of the primal problem (\ref{primaljihe4}) is linear on $\Phi(u)$. \upshape

We end this subsection with a strengthen version of Lemma \ref{projec11}.

\cor \label{projec12} Let $S\in\mathbb{R}^{r\times r}$ be a nontrivial symmetric projection matrix. For every $u\in \operatorname{int}(S\varTheta_D)$, there exists a vector $v\in S\varTheta_P$ such that $Sv\in\Phi(Su)$. An analogous result holds for $\Psi$. \upshape

\subsection{More properties}
In this subsection we give more properties of the set-valued mappings.

\lem \label{convex1}{\upshape(\cite[Corollary 3.13]{YLG20})} Both $\Phi(u)$ for every $u\in \operatorname{dom}(\Phi)$ and $\Psi(v)$ for every $v\in \operatorname{dom}(\Psi)$ are convex and closed-valued. \upshape

Two direct consequences of Lemma \ref{convex1} are given. Corollary \ref{linear1} may reference Remark \ref{remgeo1}.

\cor \label{linear1} Let $\bar{u}\in\varTheta_D$ be given. Then $\Phi(\bar{u})$ is a set if and only if the optimal objective function of the problem (\ref{primaljihe4}) is linear with respect to $v\in \Phi(\bar{u})$; An analogous result holds for $\Psi$. \upshape

\cor \label{contin1} Both $\Phi$ and $\Psi$ are outer semicontinuous. \upshape

Let $\mathscr{Y}(u)$ and $\mathscr{Y}^*(u)$ denote the sets of all feasible and optimal solutions of the problem (\ref{primaljihe3}).

\lem \label{Skman1} For every $\bar{u}\in \operatorname{int}(\varTheta_D)$, there exists a neighborhood $U$ of $\bar{u}$ such that both sets $\mathscr{X}^*(u)$ and $\mathscr{Y}^*(u)$ are nonempty and uniformly bounded, i.e., there exist compact sets $C_1,C_2$ such that \[ x^*(u)\in C_1, \quad y^*(u)\in C_2.\] \upshape

\begin{proof} Let $(x(\bar{u});w(\bar{u}),y(\bar{u}))\in\mathscr{X}\times\mathscr{Y}(\bar{u})$ be given for every $\bar{u}\in \operatorname{int}(\varTheta_D)$. If $(x(\bar{u});w(\bar{u}),y(\bar{u}))$ is strictly feasible, then there exist a triple of strictly feasible points $(x(u); w(u),y(u))\in\mathscr{X}\times\mathscr{Y}(u)$ such that $(x(u);w(u),y(u))\rightarrow (x(\bar{u});w(\bar{u}),y(\bar{u}))$ as $u\rightarrow \bar{u}$. By Theorems \ref{classical1} and \ref{slaterk1}, there exists a neighborhood $U$ of $\bar{u}$ such that for every $u\in U$, both sets $\mathscr{X}^*(u)$ and $\mathscr{Y}^*(u)$ are nonempty. Moreover, for any $x^*(u)\in \mathscr{X}^*(u)$ and $(w^*(u),y^*(u))\in\mathscr{Y}^*(u)$, one has
\[ A (x(u)-x^*(u))=0,\quad A^T(w(u)-w^*(u)) + y(u)-y^*(u)=0\] such that
\[ \langle x(u)-x^*(u), y(u)-y^*(u) \rangle=0. \] Therefore, it follows from the complementary slackness property $\langle x^*(u),y^*(u)\rangle =0$ that
\[ \langle x^*(u),y(u)\rangle +\langle x(u),y^*(u)\rangle =\langle x(u),y(u)\rangle, \] or equivalently
\[ \frac{\|x^*(u)\|}{\|x(u)\|} \cos angle(x^*(u),y(u))+ \frac{\|y^*(u)\|}{\|y(u)\|} \cos angle(x(u),y^*(u))=\cos angle(x(u),y(u)), \] where $angle(x,y)$ denotes the angle between two vectors $x$ and $y$. Therefore, for every $u\in U$, one has
\[ \frac{\|x^*(u)\|}{\|x(u)\|} \cos angle(x^*(u),y(u))\leq 1 \quad and \quad \frac{\|y^*(u)\|}{\|y(u)\|} \cos angle(x(u),y^*(u))\leq 1.\]
From the strictly feasibility of $x(u)$ and $(w(u),y(u))$, one has
 \begin{eqnarray*} c_1 &=& \max \{\cos angle(x^*(u),y(u))| (w(u),y(u))\in \mathscr{Y}(u)\}>0, \\ c_2&=& \max \{\cos angle(x(u),y^*(u))| x(u)\in \mathscr{X}\}>0. \end{eqnarray*} Then for a small enough neighborhood $U$, there exists $\varepsilon>0$ such that for some $(\hat{x}(u); \hat{w}(u),\hat{y}(u))\in\mathscr{X}\times\mathscr{Y}(u)$, one has \[c_1> \cos angle(x^*(\bar{u}),\hat{y}(\bar{u}))-\varepsilon>0 \quad \operatorname{and}\quad c_2= \cos angle(\hat{x}(\bar{u}),y^*(\bar{u}))-\varepsilon>0.\] Therefore, one has
\[ \|x^*(u)\| \leq \frac{\|\hat{x}(u)\|}{c_1}< \frac{\|\hat{x}(\bar{u})\|+\varepsilon}{\cos angle(x^*(\bar{u}),\hat{y}(\bar{u}))-\varepsilon} \] and \[ \|y^*(u)\| \leq \frac{\|\hat{y}(u)\|}{c_2}< \frac{\|\hat{y}(\bar{u})\|+\varepsilon}{\cos angle(\hat{x}(\bar{u}),y^*(\bar{u}))-\varepsilon}, \] where $(\hat{x}(u);\hat{w}(u),\hat{y}(u))\rightarrow (\hat{x}(\bar{u});\hat{w}(\bar{u}),\hat{y}(\bar{u}))\in\mathscr{X}\times\mathscr{Y}(\bar{u})$ as $u\rightarrow \bar{u}$. Therefore, $\mathscr{X}^*(u)$ and $\mathscr{Y}^*(u)$ are uniformly bounded for all sufficiently small neighborhoods $U$. \end{proof}

The idea of the proof comes from the paper \cite{YW21}. Using this result we can conclude the continuity of $\Phi$ and $\Psi$ from the uniqueness condition.

\cor \label{conu2} If $\Phi(u)$ is single-valued at $\bar{u}\in \operatorname{int}(\varTheta_D)$, then $\Phi(u)$ is continuous at $\bar{u}$ relative to $ \operatorname{int}(\varTheta_D)$. An analogous result holds for $\Psi$. \upshape

\begin{proof} The proof is from Lemma \ref{convex1}, local boundedness of $\Phi$ and $\Psi$, and \cite[Corollary 8.1]{Hoh73}. \end{proof}

For a more general result of continuity, see Theorem \ref{continu1} below.

\thm \label{xinzd2} Both $-\Phi$ and $-\Psi$ are monotone on $\mathbb{R}^r$. That is, either for any $v_i\in\Phi(u_i)$ $(i=1,2)$ or for any $u_i\in\Psi(v_i)$ $(i=1,2)$, one has
\begin{equation}\label{dandxi1} \langle u_2- u_1, v_1- v_2\rangle \geq 0. \end{equation} \upshape

\begin{proof} We only prove that $-\Phi$ is monotone. From the optimality of $x^*(u_1)$ and $x^*(u_2)$, one has
 \[ \langle c+M^Tu_1,x^*(u_1) \rangle \leq \langle c+M^Tu_1,x^*(u_2) \rangle\]
 and
 \[ \langle c+M^Tu_2,x^*(u_1) \rangle \geq \langle c+M^Tu_2,x^*(u_2) \rangle. \]
 Subtracting the first inequality from the second inequality to get
 \[\langle M^Tu_2-M^Tu_1,x^*(u_1)-x^*(u_2)\rangle \geq 0, \] or equivalently, for any $v_1\in\Phi(u_1)$ and $v_2\in\Phi(u_2)$,
 \[\langle u_2- u_1, v_1- v_2\rangle=\langle u_2- u_1,Mx^*(u_1)-Mx^*(u_2)\rangle \geq 0. \]
The proof is finished. \end{proof}

In SDP, Goldfarb and Scheinberg \cite{GS99} gave a similar result of Theorem \ref{xinzd2}. They discussed a parametric SDP problem in which the objective function depends linearly on a scalar parameter.

The following theorem plays an important role in the next analysis.

\thm \label{semid1} If $K$ is a semialgebraic set, then both the mappings $\Phi$ and $\Psi$ are semialgebraic. \upshape

\begin{proof} If $K$ is a semialgebraic set, then so is $K^*$ by Carath$\acute{e}$odory's lemma (Ref. \cite{Sch18}). From Theorem \ref{TST1}, both $\varTheta_P$ and $\varTheta_D$ are semialgebraic such that their Cartesian product $\varTheta_P\times\varTheta_D$ is semialgebraic. Finally the desired result is proved since the complementary slackness property (\ref{kktcond3}) is a polynomial equation. \end{proof}

By Theorems \ref{TST2} and \ref{semid1}, there exists a partition of the domain of $\Psi$ into semialgebraic sets $\mathcal{V}_1, \mathcal{V}_2,\cdots, \mathcal{V}_k$, so that for each $i$ and any $v_i\in \mathcal{V}_i$,
\begin{equation} \label{tongt1} \mathcal{V}_i\times (\{v_i\}\times \Psi(v_i))\cong \mathcal{V}_i\times \Psi(\mathcal{V}_i).\end{equation} It can be used as a partition of $\varTheta_P$. A detailed discussion will be given in the next section. For brevity, in the following statements we always assume that the cone $K$ is a semialgebraic set.

\section{Optimal partition} \label{parti14}
In this section we exploit the semialgebraic property of the set-valued mappings $\Phi$ and $\Psi$ (Theorem \ref{semid1}) to develop the optimal partition technique for semialgebraic mpCLOs. In Appendix B, we present an example to illustrate the behavior of the optimal partition.

\subsection{Linearity sets}
Motivated by the observation from Corollary \ref{linear1}, we review the notation of a linearity set from \cite{YLG20}.

\defn \label{djkyi1} A simply connected subset $\mathcal{V}$ of $\varTheta_P$ is called a linearity set if

 (1) $\Psi$ is not a one-to-one mapping on $\mathcal{V}$, and

 (2) $\Psi(v_1)=\Psi(v_2)$ for all $v_1,v_2\in\mathcal{V}$.

\noindent For $\varTheta_D$, the definition of the linearity set are similar and omitted. Sometimes, a linearity set of $\varTheta_P$ (or $\varTheta_D$) is called the primal (or dual) linearity set. \upshape

In the previous definition, the simply connected set $\mathcal{V}$ is always maximal. That is, for any simply connected set $G \supsetneq\mathcal{V}$, there is a point $v_0\in G-\mathcal{V}$ such that $\Phi(v_0)\ne \Psi(\mathcal{V})$. In addition, a singleton linearity set $\mathcal{V}$ is allowed.

In what follows, we need the following auxiliary results, in which some analogous results are not stated for $\varTheta_D$ or $\Phi$.

\lem \label{setjiaoh1} If $\mathcal{V}$ is a linearity set of $\varTheta_P$, then for any $u\in \Psi(\mathcal{V})$, one has $\mathcal{V}\subset \Phi(u)$ and $\operatorname{cl}(\mathcal{V})= \Phi(u)$. \upshape

\begin{proof} Immediately from Lemmas \ref{invers1} and \ref{convex1}. \end{proof}

\lem \label{setjiaoh2} Every point of the closure of linearity sets belongs to a linearity set.
\upshape

\begin{proof} Let $\mathcal{V}$ be a linearity subset of $\varTheta_P$. By Lemma \ref{setjiaoh1}, for any $u\in \Psi(\mathcal{V})$, one has $\mathcal{V}\subset \Phi(u)$ and $\operatorname{cl}(\mathcal{V})= \Phi(u)$. Let $v\in\Phi(u)$. We can assume that $v\notin\mathcal{V}$ since if not, then the result is trivial. By Lemma \ref{invers1}, there is $u_1\ne u$ such that $u_1\in \Psi(v)$ but $u_1\notin\Psi(\mathcal{V})$. By Lemma \ref{convex1}, $\operatorname{conv}\{u,u_1\}\in \Psi(v)$, i.e., $\Psi(v)$ is not a singleton set. Then from Definition \ref{djkyi1}, $v$ belongs to a linearity set that might contain only one point $v$. The same result can be used for $\varTheta_D$. \end{proof}

\lem \label{tjia1} Every point of the image of a primal (or dual) linearity set under the set-valued mapping belongs to a dual (or primal) linearity set.
\upshape

\begin{proof} Immediately from Lemmas \ref{setjiaoh1} and \ref{setjiaoh2}. \end{proof}

\lem \label{tjia2} If $\Psi(v_1)\cap\Psi(v_2)\ne \emptyset$, then there is a linearity set $\mathcal{V}$ of $\varTheta_P$ such that $v_1\in\operatorname{cl}(\mathcal{V})$ and $v_2\in \operatorname{cl}(\mathcal{V})$. \upshape

\begin{proof} If $u_0\in \Psi(v_1)\cap\Psi(v_2)$, then from Lemmas \ref{invers1} and \ref{convex1}, one has $\operatorname{conv}\{v_1,v_2\}\subset \Phi(u_0)$. By Definition \ref{djkyi1}, there is a linearity set $\mathcal{V}$ of $\varTheta_P$ such that $\Phi(u_0)\subset\operatorname{cl}(\mathcal{V})$, which implies the result. \end{proof}

\cor \label{close1} The union of all linearity sets of a conic representable set is closed. \upshape

\begin{proof} Immediately from Lemma \ref{setjiaoh2}. \end{proof}

It is easy to see that a polyhedron only contains linearity sets, or equivalently, the linearity set decomposition provides an optimal partition of a polyhedron. By lemma \ref{convex1}, a nontrivial linearity set for LPs with a scalar parameter is an invariancy interval, in which the optimal partition keeps the optimal basis decomposition unchanged. In a linearity interval, one of the optimal solutions of a pair of primal-dual LPs remains unchanged, although their objective value functions are linear. This approach identifies the range of parameters where the optimal partition remains invariant, which is consistent with the idea of Alder and Monterio \cite{AD92} and Jansen et al. \cite{JR03}.

\cor \label{bansem1} The interior and the closure of a linearity subset are semialgebraic. \upshape

\begin{proof} Let $\mathcal{V}$ be a linearity set of $\varTheta_P$. By Lemma \ref{setjiaoh1}, for any $u\in \Psi(\mathcal{V})$, one has $\mathcal{V}\subset \Phi(u)$ and $\operatorname{cl}(\mathcal{V})= \Phi(u)$. By Theorems \ref{conec1} and \ref{conec3} and Lemma \ref{convex1}, one has $\operatorname{int}(\mathcal{V}) = \operatorname{int}(\Phi(u))$ since $\mathcal{V}$ is a maximal simply connected set. By Theorems \ref{TST2} and \ref{semid1}, $\Phi(u)$ is semialgebraic such that $\operatorname{int}(\mathcal{V})$ and $\operatorname{cl}(\mathcal{V})$ are semialgebraic. \end{proof}

The following theorem provides a basic characterization of linearity sets.

\thm \label{linritys1} If $\mathcal{V}$ is a linearity subset of $\varTheta_P$, then the optimal objective function of the primal problem (\ref{primaljihe4}) is linear on $\mathcal{V}$. \upshape

\begin{proof} We may assume that $\mathcal{V}$ is not a singleton set since if it is, the result is trivial. By Definition \ref{djkyi1}, there is $\bar{u}\in\Psi(\mathcal{V})$ such that $\Phi(\bar{u})=\mathcal{V}$. Then by Corollary \ref{linear1}, the desired result holds. \end{proof}

\rem \label{linritys1pl} The inverse of Theorem \ref{linritys1} is also true. In indeed, if the optimal objective function of the primal problem (\ref{primaljihe4}) is linear on $\mathcal{V}$, then by Definition \ref{djkyi1}, Lemma \ref{tjia1} and Corollary \ref{linear1}, there is a linearity subset $\mathcal{V}'$ of $\varTheta_P$ such that $\operatorname{int}(\mathcal{V})=\operatorname{int}(\mathcal{V}')$. \upshape

\cor \label{tplog1} {\upshape \cite[Corollary 5.4]{YLG20}} If $\mathcal{V}$ is a linearity set of $\varTheta_P$, then $\mathcal{V}$ is convex, and for all $v\in \operatorname{cl}(\mathcal{V}) $ and $u\in \Psi(\mathcal{V})$, one has $v\in\Phi(u)$. \upshape

For LPs, Ghaffari Hadigheh et al. \cite{GG08} presented a similar result for the bi-parametric optimal partition.

\subsection{Nonlinearity sets}
From Theorems \ref{conec2} and \ref{conec3}, we may assume that each $\mathcal{V}_i$ in the formula (\ref{tongt1}) is simply connected. This derives us to give the following definition.

\defn \label{djkyi2} A maximal simply connected open subset of a conic representable set is called a nonlinearity set if it doesn't intersect any linearity sets. Both a linearity set and a nonlinearity set are collectively called an invariancy set. \upshape

In \cite{YLG20}, authors gave another equivalently definition and discussed the continuity of $\Phi$ and $\Psi$ on nonlinearity sets. The following corollary confirms this.

\cor \label{lianxu1} The set-valued mappings $\Phi$ and $\Psi$ are one-to-one and continuous on every nonlinearity set. \upshape

\begin{proof} By Lemmas \ref{tjia1} and \ref{tjia2}, both $\Phi$ and $\Psi$ are one-to-one on every nonlinearity set. Then by Corollary \ref{conu2}, they are also continuous. \end{proof}

The following corollary shows that a nonlinearity set with a scalar parameter is an invariancy interval. So Definition \ref{djkyi2} extends the concept given by Mohammad-Nezhad and Terlaky \cite{MT21} in the case of SDP.

\cor \label{open1} \cite[Theorem 5.2]{YLG20} Every nonlinearity set is open. \upshape

\cor \label{bansem2} The interior and the closure of a nonlinearity subset are semialgebraic. \upshape

\begin{proof} Let $\mathcal{U}$ denote the union of all nonlinearity sets of $\varTheta_D$. By Corollary \ref{lianxu1}, one has \[\mathcal{U}=\{u\in\mathbb{R}^r| \operatorname{ the\ system\ (\ref{kktcond10})-(\ref{kktcond30})\ has \ a \ unique \ solution} \ (\bar{x};\bar{w},\bar{y}) \}. \] Clearly, $\mathcal{U}$ is the projection of a Boolean combination of polynomial equalities and inequalities, and by Theorem \ref{TST1}, $\mathcal{U}$ is a semialgebraic set. Then by Corollaries \ref{tjia1} and \ref{lianxu1}, $\Phi(\mathcal{U})$ is also semialgebraic. This concludes the final result. \end{proof}

The main result of this paper is as follows.

\thm \label{bansem3} A conic representable set is a union of finite linearity and/or nonlinearity sets. \upshape

\begin{proof} By Theorems \ref{semid1} and \ref{TST2} or \ref{conec2}, the formula (\ref{tongt1}) holds such that $(\mathcal{V}_1,\mathcal{V}_2,\cdots,\mathcal{V}_k)$ is a partition of $\varTheta_P$. Then from Corollaries \ref{bansem1} and \ref{bansem2}, every $\mathcal{V}_i(i=1,2,\cdots,k)$ is either a linearity set or a nonlinearity set. The proof is completed. \end{proof}

The crucial component of the above proof is the observation that every invariancy set satisfies the formula (\ref{tongt1}). We say that such an optimal partition is the finest. That is, if a partition of $\varTheta_P$ into the invariancy sets $\mathcal{V}_1, \mathcal{V}_2,\cdots, \mathcal{V}_k$ is defined as Definitions \ref{djkyi1} and \ref{djkyi2}, then there is not a series of the disjoint semialgebraic sets $\mathcal{V}'_1, \mathcal{V}'_2,\cdots, \mathcal{V}'_l$, so that the following three conditions are true:

(1) $\operatorname{cl}\left(\bigcup\limits_{j=1}^l\mathcal{V}_j'\right)\supset\varTheta_P$;

(2) For each $j$ and any $v_j\in \mathcal{V}_j'$,
\[\mathcal{V}_j'\times (\{v_j\}\times \Psi(v_j))\cong \mathcal{V}_j'\times \Psi(\mathcal{V}_j');\]

(3) There are two indices $i_0\in\{1,2,\cdots,k\}$ and $j_0\in\{1,2,\cdots,l\}$ such that $\mathcal{V}_{j_0}'$ is a proper subset of $\mathcal{V}_{i_0}$.

\cor \label{touying1} If a conic representable set contains a nonlinearity set, then its nontrivial projection contains a nonlinearity set. \upshape

\begin{proof} This is an easy consequence of Corollary \ref{projec12}. \end{proof}

Just like Theorem \ref{linritys1}, the following result gives a basic characterization of nonlinearity sets.

\thm \label{linritys2} If $\mathcal{V}$ is a nonlinearity set of $\varTheta_P$, then the optimal objective function of the primal problem (\ref{primaljihe4}) is nonlinear on $\mathcal{V}$. \upshape

\begin{proof} This result follows from Theorem \ref{linritys1}, Remark \ref{linritys1pl} and Definition \ref{djkyi2}. \end{proof}

In Example \ref{patex10}, $\varTheta_D$ contains only one nonlinearity interval $(0,+\infty)$, in which the optimal objective function $2\sqrt{u}$ is nonlinear and smooth.

Theorems \ref{linritys1} and \ref{linritys2} show that there are a clear distinction between the linearity and nonlinearity sets. Unlike a linearity set, for every point $u$ of a nonlinearity set of $\varTheta_D$, both the supporting hyperplane $H_u$ defined by (\ref{zcsup1}) and the cutting set $C_v$ defined by (\ref{zcsup2}) are unique. In other words, only at the unique extreme point $x^*(u)$, the unique supporting hyperplane $H_u$ is tangent to $\mathscr{X}$, whereas the unique cutting set $C_v$ intersects $\mathscr{X}$. Therefore, the feasible set $\mathscr{X}$ is smooth at the extreme point $x^*(u)$.

\cor \label{ghui1} Let $a\in\mathbb{R}^q$ denote the first column vector of the transpose of the matrix $A$, and $b_1\in \mathbb{R}$ denote the first component of the vector $b$. If the boundary of the set $\mathscr{X}_1=\{x\in K| a^Tx=b_1\}$ is smooth, then the nonlinearity sets of $\varTheta_P$ and $\varTheta_D$, defined by (\ref{biaoshi10}) and (\ref{biaoshi20}), respectively, are not empty. \upshape

\begin{proof} Let the matrix $A$ be decomposed as follows
\[ A = \left(\begin{array}{c} a^T \\ A_0 \end{array} \right). \] Consider the following problem
\begin{equation} \label{prim10c} \begin{array}{ll} \min & \left\langle c+ \left( A_0^T, M^T \right) \left(\begin{array}{c} u_0 \\ u \end{array} \right),x \right\rangle \\ s.t. & a^Tx=b_1,\\ & x\in K \end{array} \end{equation} and the corresponding conic representable set \[\varTheta_P^0=\left\{ \left. \left(\begin{array}{c} v_0 \\ v \end{array} \right)\right|\operatorname{exists} \ \bar{x}\ \operatorname{s.t.} \ a^T \bar{x}=b_1, \left(\begin{array}{c} A_0 \\ M \end{array} \right) \bar{x}= \left(\begin{array}{c} A_0 \\ M \end{array} \right)d + \left(\begin{array}{c} v_0 \\ v \end{array} \right), \bar{x}\in K \right\},\] where $d\in\mathbb{R}^q$ is feasible for the problem (\ref{prim10c}). If $A_0d=(b_2,\cdots,b_m)^T$ and \[S= \left(\begin{array}{cc} 0_{m-1} & 0 \\ 0& I_r \end{array} \right),\] then the problem (\ref{primaljihe1}) is a projection of the problem (\ref{prim10c}).
By hypothesis, the nonlinearity set of $\varTheta_P^0$ exists. By Corollary \ref{touying1}, the nonlinearity set of $\varTheta_P=S\varTheta_P^0$ exists. Finally, by Corollary \ref{lianxu1}, the nonlinearity set of $\varTheta_D$ also exists. \end{proof}

Corollary \ref{ghui1} shows that a spectrahedon is usually smooth, which answers the open problem proposed by Hauenstein et al. \cite{Hau19}.

\subsection{Transition faces}
In this subsection we focus on the set that separate invariancy sets. It can be defined as follows.

\defn \label{djkyi3} An invariancy set of a conic representable set is called a transition face if it is not full-dimensional in $\mathbb{R}^r$. In particular, if its affine dimension is zero, then the transition face reduces a transition point; and if its affine dimension is equal to 1, then the transition face reduces a transition line, and etc. \upshape

\cor \label{trans2} Every transition face is a linearity set. \upshape

\begin{proof} By Corollaries \ref{close1} and \ref{open1}, the boundary of a invariant set does not intersect any nonlinearity sets, which yields the result. \end{proof}

We refer to a linearity set as a nontrivial linearity set if it is not a transition face. In contrast, a transition face is called a trivial linearity set.

\cor \label{nonlinear2} For every boundary point of a nonlinearity set, either it is an element of a transition face, or the set-valued mapping of that point is undefined. \upshape

\begin{proof} If the set-valued mapping on the boundary point of a nonlinearity set is well defined, then by Corollary \ref{trans2}, the point belongs to a transition face. \end{proof}

\thm \label{continu1} Both $\Phi$ and $\Psi$ are continuous on $\operatorname{int}(\varTheta_D)$ and $\operatorname{int}(\varTheta_P)$, respectively. \upshape

\begin{proof} By Corollaries \ref{contin1} and \ref{lianxu1}, we only need prove that both $\Phi$ and $\Psi$ are inner semicontinuous on every linearity set. If the linearity set intersects a nonlinearity set, then this result follows from Corollary \ref{nonlinear2}. Otherwise, it follows from Corollary \ref{trans2}. \end{proof}

Now we focus on vertices of the feasible sets. Such vertices are naturally among the first object to understand in a study of the boundary structure.

\cor \label{trans3} The invariant set associated with a vertex is a transition face. \upshape

\begin{proof} For every fixed $\bar{u}\in\varTheta_D$, if the corresponding extreme point $x^*(\bar{u})$ is a vertex, then $\Phi(\bar{u})$ is well defined, and the normal cone of $\mathscr{X}$ at $x^*(\bar{u})$ is full-dimensional. Therefore, by Lemma \ref{tjia2}, $\bar{u}$ belongs to a transition face of $\varTheta_D$. Similar arguments are applied to the extreme points of $\mathscr{Y}$. \end{proof}

For a single parameter optimization, i.e., $r=1$, both the nontrivial linearity set and the nonlinearity set are an interval, and the transition face reduces a point.

\cor \label{translem1} If $r=1$, then $\varTheta_P$ contains finitely transition points, and the invariant set associated with a vertex of $\mathscr{X}$ is a transition point. \upshape

\begin{proof} The former is a direct consequence of Theorem \ref{bansem3}, and the latter is from Corollaries \ref{trans2} and \ref{trans3}. \end{proof}

A similar result for the second-order conic optimization was given by Mohammad-Nezhad and Terlaky \cite{MT21}. In LP, every transition point is associated with only one vertex of the polyhedron. However, in SDP, a transition point is not necessarily associated with a vertex of the spectrahedron. A counterexample is given in Appendix B.

\thm \label{trans5} The vertices of the feasible set $\mathscr{X}$ are finite. \upshape

\begin{proof} Assume otherwise. Consider the problem (\ref{primaljihe4}) with $r=q-m$. By Corollary \ref{trans3}, the set $\varTheta_P$ covers the transition points associated with all vertices of the feasible set $\mathscr{X}$. By Theorem \ref{bansem3}, there exists a partition of the domain of $\varTheta_P$ into invariancy sets $\mathcal{V}_1, \mathcal{V}_2,\cdots,\mathcal{V}_{k}$, in which each is semi-algebraic. Then there is an index $i_0\in\{1,2,\cdots,k\}$ such that $\mathcal{V}_{i_0}$ is associated with an infinite number of vertices of $\mathscr{X}$. Furthermore, there is a symmetric projective matrix $S\in\mathbb{R}^{r\times r}$ with rank 1 such that $S\mathcal{V}_{i_0}$ is associated with an infinite number of vertices of $\mathscr{X}$. By Lemma \ref{projec12}, there exists an infinite number of the transition points in $S\mathcal{V}_{i_0}$, which results in a contradiction with the last remark. \end{proof}

\section{Conclusions}
This paper deals with a system of conic linear equations which is called the mpKKT property and provide some characterizations of its solution set. Such a system is an extension of KKT systems derived from LP and SDP into the case with semialgebraic convex cone $K$. This extends the work presented in a previous draft \cite{YLG20}, where we first explored a novel notion of duality in mpCLOs, defining and using conic representable sets and set-valued mappings defined on these. We then introduce the notions of linearity and nonlinearity sets of the optimal partition in the general case of mpCLOs and show that a conic representable set is a union of finite linearity and/or nonlinearity sets. This generalizes the work presented in \cite{Hau19,MT20}, where the notion of a nonlinearity interval was first introduced in the case of SDPs. Finally, we investigate the boundary structure of the feasible set of semialgebraic conic linear optimization problems, proving that the number of vertices of the feasible set is finite.

We establish more useful properties for the mpCLOs in which the perturbation of the objective function occurs in many fixed directions. Many of the results of this paper are heavily dependent on the Slater type condition. When this condition vanishes, the set-valued mappings defined by (\ref{ghjj1})-(\ref{ghjj2}) become much more complicated. We leave this case for future investigation. Another interesting issue of future investigation is to provide more concrete properties on the boundary structure. The work in \cite{MST16} also provides some interesting results in this direction.

\appendix
\section{Proof of Theorem \ref{slaterk1}}
In this appendix, we prove Theorem \ref{slaterk1}. We start with several auxiliary lemmas.

\lem \label{dua1} (\cite{NN94}) Let $D$ be a closed cone in $\mathbb{R}^q$, and $D^*$ be its dual cone. Then

(1) $D$ is point if and only if $D^*$ has a nonempty interior. $M$ is a closed pointed cone with a nonempty interior if and only if the dual cone has the same properties.

(2) $a\in \operatorname{int}(D^*)$ if and only if $\langle a,x\rangle>0$ for all nonzero vectors $x\in D$. \upshape

\lem \label{dua2} Let $c,a\in \mathbb{R}^q$ and $\{\mu_k\}_{k=1}^{\infty}\subset I_{>}$ be given, where
 \[ I_{>} =\{\lambda| c+\lambda a\in \operatorname{int}(K^*) \}. \] If $\lim\limits_{k\rightarrow +\infty}\mu_k=0$, then $c\in K^*$. \upshape

\begin{proof} From Lemma \ref{dua1}, $K^*$ is a pointed, closed, convex, solid cone such that $(K^{*})^{*}=K$; and for any nonzero vectors $x\in K$, one has
\[\langle c+\mu_ka,x\rangle >0, \qquad k=1,2,\cdots. \] Letting $k$ tend to infinity, one get \[\langle c,x\rangle \geq 0. \] Then $c\in K^*$ as desired. \end{proof}

\lem \label{openint1} Let $c,a\in \mathbb{R}^q$ be given. Then the set $I_{>}$ is an interval. \upshape

\begin{proof} We show that $I_{>}$ is convex and thus an interval. Assume that $\lambda_1$ and $\lambda_2$ belong to $I_{>}$ and that $\lambda_1<\lambda_2$. Then for any $\theta\in (\lambda_1,\lambda_2)$,
\[ c+\theta a= \frac{\lambda_2-\theta}{\lambda_2-\lambda_1} (c+\lambda_1 a)+\frac{\theta-\lambda_1}{\lambda_2-\lambda_1} (c+\lambda_2 a) \in \operatorname{int}(K^*).\] Thus, $\theta\in I_{>}$ as desired. \end{proof}

If $c\in \operatorname{int}(K^*)$, then $I_{>}$ is a nonempty open interval. $I_{>}$ is a finite interval $(\lambda_l,\lambda_u)$ if $\pm a\notin K^*$. Indeed, if there a big enough real number $\lambda>0$ such that $c+\lambda a\in \operatorname{int}(K^*)$, then by Lemma \ref{dua2}, $\frac{1}{\lambda}c+a \in \operatorname{int}(K^*)$ implies that $a\in K^*$, which results in one conflict. Or if there a small enough real number $\lambda<0$ such that $c+\lambda a\in \operatorname{int}(K^*)$, then $-\frac{1}{\lambda}c-a \in \operatorname{int}(K^*)$ results in another conflict. In a similar manner we can prove that $I_{>}$ is of the form $(\lambda_l,+\infty)$ if $a\in \operatorname{int}(K^*)$, while $I_{>}=(-\infty,\lambda_u)$ if $-a\in \operatorname{int}(K^*)$. Here $c+\lambda_l a\in (K^*-\operatorname{int}(K^*))$ and/or $c+\lambda_u a\in (K^*-\operatorname{int}(K^*))$. This yields the following result.

\lem \label{openint2} Let $c,a\in \mathbb{R}^q$ be given. If $c\in \operatorname{int}(K^*)$, then there is a real number $\lambda\in\mathbb{R}$ such that $c+\lambda a\in (K^* -\operatorname{int}(K^*))$. \upshape

Let $S\in\mathbb{R}^{r\times r}$ be a symmetric projection matrix. If $M$ is replaced by $SM$, then the mpKKT property (\ref{kktcond1})-(\ref{kktcond3}) becomes
\begin{subequations}
\begin{align}\label{skktcond1} &A x =b,\quad SMx= SMd + v,\quad x\in K, \\ \label{skktcond2} &
A^Tw+y=c+M^TSu,\quad w\in\mathbb{R}^m, \quad y\in K^*,\\ \label{skktcond3}& \langle x,y \rangle=0.
\end{align}
\end{subequations} Clearly, it is the KKT property of a pair of primal-dual problems
 \begin{equation} \label{primaljihe40c} \begin{array}{ll} \min & \langle c, x \rangle \\ s.t. & Ax=b,\\ & SMx=SMd+ v, \\ & x \in K \end{array} \end{equation} and
\begin{equation} \label{primaljihe20c} \begin{array}{ll} \max & b^Tw \\ s.t. & A^Tw+y=c+M^TSu, \\ & w\in\mathbb{R}^m, y \in K^*. \end{array} \end{equation}
Therefore, the problems (\ref{primaljihe40c}) and (\ref{primaljihe20c}) denote the projections of the problems (\ref{primaljihe4}) and (\ref{primaljihe3}), respectively.

For the proof of Theorem \ref{slaterk1} we need the following key lemma.

\lem \label{projec11} Let $S\in\mathbb{R}^{r\times r}$ be a nontrivial symmetric projection matrix. Suppose that $\Phi(u)$ is well defined for some $u\in \varTheta_D$. If $u\in S\varTheta_D$, then for any $v\in \Phi(u)$, one has $Sv=v$. An analogous result holds for $\Psi$. \upshape

\begin{proof} Firstly, it follows from $u\in S\varTheta_D$ that $Su=u$. If $v\in \Phi(u)$, then the system (\ref{skktcond1})-(\ref{skktcond3}) holds. Leftmultipling both sides of the second equality of (\ref{skktcond1}) by $S$, one has
\[ SMx= SMd + Sv. \] Subtract the two equations and get $Sv=v$. The desired result for $\Phi$ is proved.

Secondly, note that the problem
\begin{equation} \label{prijihe20c} \begin{array}{ll} \max & b^Tw+(Md+Sv)^Ts \\ s.t. & A^Tw+M^Ts+y=c, \\ & w\in\mathbb{R}^m,s\in\mathbb{R}^r, y \in K^* \end{array} \end{equation} can be rewritten as \begin{equation} \label{prijihe20c} \begin{array}{ll} \max & b^Tw+(SMd+v)^T(Ss) \\ s.t. & A^Tw+(SM)^T(Ss)+y=c, \\ & w\in\mathbb{R}^m,s\in\mathbb{R}^r, y \in K^*. \end{array} \end{equation} If $v\in S\varTheta_P$, then $Sv=v$ such that $Ss=s$. Therefore, one has $Su=u$. The proof is completed. \end{proof}

{\it Proof of Theorem \ref{slaterk1}.} Let us prove the result by the mathematical induction for $m$.

{\bf Initial step:} By Lemma \ref{openint2}, this result holds for $m=1$.

{\bf Inductive step:} By Remark \ref{strong1}, we can assume that the mpKKT property (\ref{kktcond1})-(\ref{kktcond3}) holds for all $1\leq m\leq k$ and $r\leq q-m$, in which $u\in \operatorname{int}(\varTheta_D)$, and the rank of the matrix $A$ is less than $k$. We may choose $m$ and $r\geq 2$ such that $m+r-\operatorname{rank}(S)=k+1\leq q$, where $S\in\mathbb{R}^{r\times r}$ is a symmetric projection matrix such that $Su=u$. By Lemma \ref{projec11}, one has $Sv=v$. Therefore, the mpKKT property (\ref{kktcond1})-(\ref{kktcond3}) holds for \[ A:= \left( \begin{array}{c} A \\ (I_r-S)M \end{array} \right), \quad b: = \left( \begin{array}{c} b \\ (I_r-S)Md \end{array} \right) \quad \operatorname{and} \quad M: = SM, \] in which the rank of the updated matrix $A$ increases to $k+1$. This concludes the result by induction. $\square$

Note that the strong duality theorem of conic programming is not used in the above proof.

\section{The behavior of optimal partition}
In this section we give an example to illustrate the optimal partition of conic representable sets. Our recent paper \cite{YLG20} offers the following example, in which the complex calculations are omitted.

\exam \label{patex2} Consider the bi-parameter SDP problem
\begin{equation} \label{mostl1} \begin{array}{ll} \min\limits_{x\in\mathbb{S}^3_+} & -x_{23}+ u_1 x_{12} -u_2x_{13} \\ s.t. & x_{ii} =1,\quad i=1,2,3, \end{array} \end{equation} in which $\varTheta_D=\mathbb{R}^2$ and \[\varTheta_P=\{(v_1,v_2)^T| -1\leq v_1,v_2\leq 1\}.\] \upshape

The following parallel table lists the optimal partition and the values of the set-valued mappings $\Phi(u)$ and $\Psi(v)$.
\[\begin{array}{c} Table \ 1 \quad The \ set-valued\ mappings \ and \ the \ optimal\ partition \ in \ Example \ \ref{patex2} \\ \begin{array}{||c||c|c|c|c|c||c||c||} \hline \hline v & \Psi(v) \\ \hline (1,1)^T & \varTheta_D^1=\{(u_1,u_2)^T |u_1< -1, u_2< -1, u_1+u_2+u_1u_2\geq0 \} \\ \hline(-1,-1)^T & \varTheta_D^2=\{(u_1,u_2)^T|u_1> 1, u_2> 1, -u_1-u_2+u_1u_2\geq0 \} \\ \hline (1,-1)^T & \varTheta_D^3=\{(u_1,u_2)^T|u_1< 1, u_2> -1, u_2\geq u_1, u_2-u_1-u_1u_2\geq0 \} \\ \hline\varTheta_P^l &(0,0)^T=\varTheta_D^3\cap\varTheta_D^4 \\ \hline (-1,1)^T&\varTheta_D^4 = \{(u_1,u_2)^T|u_1> -1, u_2< 1, u_1\geq u_2, u_1-u_2-u_1u_2\geq0 \} \\ \hline\varTheta_P^0 & \varTheta_D^0= \mathbb{R}^2 -\bigcup\limits_{i=1}^4\varTheta_D^i \\ \hline \Phi(u) & u \\ \hline \hline\end{array} \end{array}\] in which $\varTheta_P^l=\{(v_1,-v_1)^T| -1< v_1< 1\} $ and $\varTheta_P^0=\{(v_1,-v_2)^T| -1< v_1\ne v_2 < 1\}$. The set-valued mappings between the primal and dual nonlinearity sets can be expressed as \[\Phi(u)= \left( \begin{array} {ccc} \frac{u_2}{2u_1^2} - \frac{u_2}{2} -\frac{1}{2u_2} \\ -\frac{1}{2u_1} -\frac{u_1}{2} +\frac{u_1}{2u_2^2} \end{array} \right), \quad \Psi(v)= \left(\begin{array}{c}-v_2-v_1\sqrt{\frac{1-v_2^2}{1-v_1^2}} \\ -v_1-v_2 \sqrt{\frac{1-v_1^2}{1-v_2^2}} \end{array} \right). \]

\begin{figure}[htbp]
\centering\includegraphics[width=3.5in]{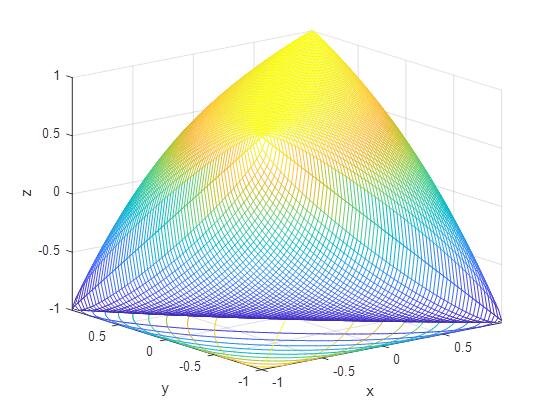}
\caption{The primal feasible set is a 3-elliptope}\label{fig:21}
\end{figure}

\begin{figure} \centering \includegraphics[height=4cm,width=5cm]{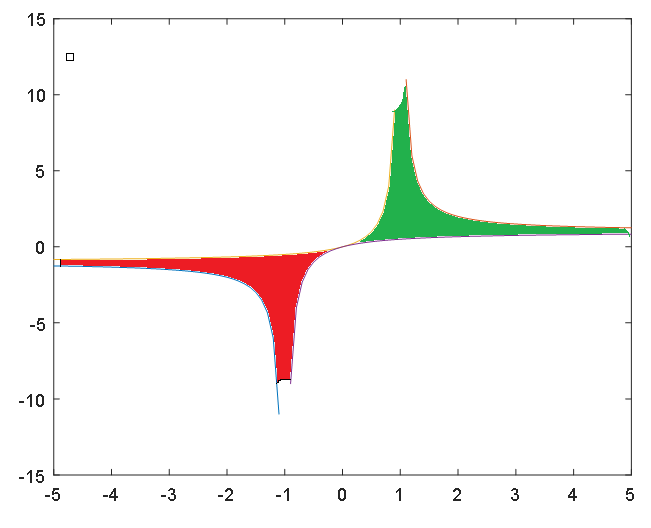} \caption{The invariancy sets of $\varTheta_D$ separated by four curves} \label{fig:3} \end{figure}

The primal conic representable set $\varTheta_P$ contains four transition points $(1,1)^T$, $(-1,-1)^T$, $(1,-1)^T$, $(-1,1)^T$, a transition line $\varTheta_P^l$, and two nonlinearity sets contained in the set $\varTheta_P^0$, in which four transition points are associated with four vertices of a 3-elliptope (see Figure \ref{fig:21}). And the dual conic representable set $\varTheta_D$ contains four nontrivial linear invariancy sets $\varTheta_D^1,\varTheta_D^2,\varTheta_D^3-\{(0,0)^T\},\varTheta_D^4-\{(0,0)^T\}$, a transition point $(0,0)^T$, two nonlinearity sets contained in the set $\varTheta_D^0$ (see Figure \ref{fig:3}).

This example shows different types of the optimal partition. In particular, if $u_1=u_2$, then the problem (\ref{mostl1}) reduces a single parameter SDP. A detail investigation can be found in \cite{MT20}. The interested readers may use the example to verify the conclusions of Section \ref{parti14}.


\bibliographystyle{amsplain}

$$$$
\end{document}